\documentstyle{article}

\def\H{{\bf H}}
\def\C{{\bf C}}
\def\R{{\bf R}}

\def\Z{{\bf Z}}

\def\Q{{\bf Q}}

\def\slz{{\rm SL}_2(\Z)}

\def\qed{\hfill\framebox[3mm][t1]{\phantom{\small x}}}

\def\VB{V\!B^{\natural}}
\def\VF{V_{\rm Fermi}}
\def\V#1{V_{(#1)}}

\def\halb{\hbox{$\frac{1}{2}$}}

\title{Self-dual Vertex Operator Superalgebras with Shadows of large minimal
weight}

\author{Gerald H\"ohn}
\date{15.~April 1997} 

\begin{document}

\bibliographystyle{../amsalpha} 

\newtheorem{satz}{Theorem}
\newtheorem{lemma}{Lemma}
\newtheorem{korollar}[satz]{Corollary}
\newtheorem{definition}[satz]{Definition}
\newtheorem{vermutung}[satz]{Conjecture}

\renewcommand{\baselinestretch}{1.2}

\maketitle

\begin{abstract} 
The shadow $V'$ of a self-dual vertex operator superalgebra $V$ is defined 
as the direct sum of the irreducible modules of its even vertex operator 
subalgebra $V_{(0)}$ not contained in $V=\V0\oplus \V1$. 
We describe the self-dual ``very nice'' unitary rational
vertex operator superalgebras $V$ of rank $c$ whose shadows have
the largest possible minimal weights $\frac{c}{8}$ or $\frac{c}{8}-1$.
The results are analogous to and imply the corresponding results for self-dual
binary codes and lattices.
\end{abstract}

\section{Introduction}

There exist a deep not completely understood analogy between codes, lattices and 
vertex operator superalgebras, as studied for example in~\cite{Go-mero,Ho-dr}.
In this short note we establish a further example of this analogy.
Elkies has recently been given a description of self-dual binary codes and 
lattices with long shadows~\cite{El-Z,El-shadow}. 
We obtain such a description for vertex operator superalgebras (SVOAs) 
which form a mathematical algebraic definition of chiral algebras in conformal 
field theory. The miracle
is, that we are working with different objects and different numbers and 
getting results of the same structure.\footnote{The 
reader will find that the structure of the results and proofs of this
paper are exactly the same as in~\cite{El-shadow}. This is intended to emphasise
the analogy.} 
By using the constructions that attach a lattice to a code and a SVOA to a 
lattice the theorems of Elkies for codes and lattices follow from
our theorems for SVOAs. 

\smallskip
Let $V$ be a self-dual unitary ``very nice'' rational vertex operator 
superalgebra of rank $c\in\halb\Z_{\geq 0}$ defined over the complex numbers. 
(A short definition is given at the end of this introduction.
For more details see~\cite{Ho-dr}, Chapter~1 and 2.)
The even vertex operator subalgebra $\V0$ of $V=\V0\oplus\V1$ consists of the 
vectors of integral weight and has besides $\V0$ and $\V1$ one (if $c\in 
\Z+\halb$) or two (if $c\in\Z$) more irreducible modules $\V2$ (and~$\V3$). We 
call $V'=\V2$ ($c\in\Z+\halb$) resp.~$V'=\V2\oplus\V3$ ($c\in\Z$) the 
{\it shadow\/} of $V$; generalising the corresponding definitions for codes 
and lattices~\cite{El-shadow}. In the special case that $V$ is a vertex 
operator algebra (VOA), i.e.~$V=\V0$, we set $V'=V$.

In~\cite{Ho-dr}, Chapter~5 we have given, in analogy to the known results for codes 
and lattices, a description of all extremal SVOAs, i.e.~SVOAs $V$ with the 
largest possible minimal weight $\halb\left[\frac{c}{8}\right]+\halb$. 
An in some sense dual
problem is the classification of SVOAs with shadows of large minimal weight
$h(V')\in \Z+\frac{c}{8}$, this is the smallest conformal weight of a nonzero 
element in $V'$. If $V$ is the tensor product $\VF^{\otimes 2c}$ of $2c$ copies 
of the ``single fermion'' SVOA $\VF$ of rank $\halb$ one has $h(V')=
\frac{c}{8}$. We show that if $V\not\cong\VF^{\otimes 2c}$ one has $h(V')<
\frac{c}{8}$ and give a complete description of all $V$ with $h(V')=\frac{c}{8}
-1$ (up to uniqueness questions and some work in the classification of all
self-dual SVOAs of rank smaller then $24$ which is in 
preparation~\cite{Hoehn2}). For the rank $c=23\halb$, one gets
the {\it shorter moonshine module\/} $V\!B^{\natural}$.
We obtain also a lower bound on the 
dimension of the weight $1$ part $V_1$ among the $V$ with minimal weight at
least $1$.

\medskip
We give very briefly the definition of a self-dual unitary ``very nice'' 
rational SVOA (cf.~\cite{Ho-dr}, Chapter~1 and~2). The more combinatorial
results and proofs in the following sections should be understandable without 
further knowledge of VOA-theory as is needed for the discussion at the end of 
this introduction.

A {\it ``nice''\/} SVOA $V=(V,{\bf 1},\omega,Y)$ is a {\it simple\/} SVOA
$$V=\bigoplus_{n\in \halb\Z_{\geq 0}}V_n=V_{(0)}\oplus V_{(1)}, 
\quad\hbox{$V_0=\C\cdot{\bf 1}$},$$
such that $V$ and its modules are direct sums of highest weight modules for the 
Virasoro algebra
generated by the coefficients of $Y(\omega,z)$ and $V$ satisfies Zhu's $C_2$
finiteness condition ${\rm dim}(V/V_{-2}V)<\infty$. In addition to~\cite{Ho-dr}
we assume also that $V$ satisfies the conditions to use the
tensor product theory for $V$-modules (cf.~\cite{Huang}). 
Most of the important examples of VOAs are ``nice.''
{\it Rational\/} means that $V$ has only finitely many irreducible modules
and every module is complete reducible. Furthermore all the fusion rules, these 
are the dimensions of the intertwiner spaces, are assumed to be finite.
The SVOA $V$ is called {\it unitary\/} if there exists a real form $V_{\R}$ of
$V$ such that the canonical invariant bilinear form on $V_{\R}$  
is positive definite. Additionally the conformal weights of the irreducible 
modules must be positive (this property probably follows).
In the case that the odd part $V_{(1)}$ is nonzero the above properties should
respectively also hold for the even sub-VOA $V_{(0)}=\bigoplus_{n\in \Z_{\geq 0}}V_n$. 
Finally we denote $V$ as {\it self-dual\/} if the only 
irreducible $V$-module is $V$ itself. 

The rank $c$ of a self-dual ``nice'' rational SVOA is a integral 
or half-integral number (see~\cite{Ho-dr}, Th.~2.2.2). 

For a self-dual SVOA $V$ to be {\it ``very  nice''\/} the following additional
properties for its even sub-VOA $V_{(0)}$ are assumed:
\newline\phantom{123} 1.~One has $V=V_{(0)}$ or
the fusion algebra of $V_{(0)}$ is given by 
$\Z[\Z_2\times\Z_2]$ if $c\in 2\Z$, $\Z[\Z_4]$ if $c\in 2\Z+1$ and the Ising 
model fusion rules if $c\in\Z+\halb$. \newline
\phantom{123} 2.~The ~$\slz=\langle S,T \rangle$-representation on the genus 
one correlation functions is given by the diagonal matrix 
$\widetilde{T}=e^{-2\pi i \frac{c}{24}}
(1,-1,e^{2 \pi i\frac{c}{8}}, e^{2 \pi i\frac{c}{8}})$ (if $c$ is integral,
$V\not=V_{(0)}$) or $\widetilde{T}=e^{-2\pi i \frac{c}{24}}
(1,-1,e^{2 \pi i\frac{c}{8}} )$ (if $c$ is half-integral). For the values
of $\widetilde{S}$ compare~\cite{Ho-dr}, (2.7).  \newline
One more property is needed:
The tensor product theory defines on the direct sum of the irreducible 
$V_{(0)}$-modules the structure of an {\it intertwiner algebra\/} 
(see~\cite{Huang}, Th.~3.5) and this structure induces {\it fusing\/} and
{\it braiding isomorphisms\/} between tensor products of intertwiner spaces
(see~\cite{Huang}, 
Remark~3.3). In~\cite{MS}, Appendix~D 1.~and~E it was discussed that in the 
cases of the fusion algebra as above the fusing and braiding maps are
well defined and depend only on the residue class $c\ ({\rm mod}\ 8)$. Since 
we don't want to check that all the properties of the braiding and fusion 
maps as assumed in~\cite{MS} are satisfied, we simply {\it assume\/} that
\newline\phantom{123} 3.~the braiding and fusing is the one described in~\cite{MS}. 

The importance of the ``braided fusion algebra''
(this is the fusion algebra together with the braiding  and fusion isomorphisms 
and the $\slz$-representation) of a VOA $V$ is the fact
that it describes completely the possible (S)VOA-extensions $W$ of $V$:
Indeed, a vertex operator defining a (S)VOA-structure on the $V$-module $W$ 
is given by elements in the intertwiner spaces for $V$ and 
the Jacobi identity for $W$ is equivalent to 
the commutativity which can be expressed in terms of the braiding maps. 
Furthermore, the ``braided fusion algebra'' of $W$ can be described completely 
in terms of the ``braided fusion algebra'' of $V$ and the vertex operator
of $W$ as elements inside.\footnote{ The ``braided fusion algebra'' should be
equivalent to a braided tensor category with extra properties,
called ``modular category'' in~\cite{Turaev}.
This means one can
associate to a VOA a $3$-dimensional topological quantum field theory, which should
be thought as the VOA-analog of the discriminant form of an integral lattice.  
}
As an application, one sees that the SVOA-structure of the 
$V_{(0)}$-module $V=V_{(0)}\oplus V_{(1)}$ is unique up to 
isomorphism.

\section{Classification results}

The subalgebra generated by the weight-$\frac{1}{2}$-part $V_{1/2}$ of a
unitary SVOA $V=\bigoplus_{n\in\frac{1}{2}\Z_{\geq 0}} V_n$ is isomorphic to 
the tensor product $\VF^{\otimes k}$ of $k={\rm dim}\,V_{1/2}$ copies of
the SVOA $\VF$, the unique unitary SVOA with the smallest possible rank 
$\frac{1}{2}$. 
The graded vector space $\VF^{\otimes k}$ is in fact an irreducible module for 
the infinite Clifford algebra generated by the coefficients of the vertex 
operators $Y(a,z)=\sum_{n \in \bf{\scriptstyle Z}}a_n\,z^{-n-1}$ of elements 
$a\in V_{1/2}$. 
\begin{satz}\label{thVF}
For the minimal weight of the shadow of a self-dual unitary
``very nice'' rational vertex operator superalgebra $V$ of rank $c$ one
has $h(V')\leq\frac{c}{8}$ with equality exactly if 
\hbox{$V\cong\VF^{\otimes 2c}$}.
\end{satz}
{\bf Proof:} The character $\chi_M$ of a module $M$ of a vertex operator
(super)algebra of rank $c$ is defined as the generating series
$$ \chi_M=q^{-\frac{c}{24}}\sum_{n\in \Q}{\rm dim\,}M_n\cdot q^n$$ of the
dimensions of its homogeneous pieces. Taking $q=e^{2\pi i\tau}$ the function
$\chi_M(\tau)$ is for sufficient regular (S)VOAs
holomorphic in the upper half-plane of the complex numbers~\cite{zhu-dr}.
For a SVOA as in the theorem we have proved in \cite{Ho-dr} (Theorem 2.2.3):

The character $\chi_V$ is a weighted-homogeneous polynomial $P_V(x,y)$ of
weight $c$ in the characters $x=\chi_{1/2}:=\chi_{\VF}=
\sqrt{\frac{\Theta_{\Z}}{\eta}}$
(weight $1/2$) of the SVOA $\VF$ and $y=\chi_8:=\chi_{V_{E_8}}=
\frac{\Theta_{E_8}}{\eta^8}$ (weight $8$) of the lattice-VOA $V_{E_8}$. 

Here $\Theta_{\Z}$ resp.~$\Theta_{E_8}$ denote the theta series
of the lattices $\Z$  resp.~$E_8$, the root lattice of the Lie group $E_8$, and
$\eta=q^{1/24}\,\prod_{n=1}^{\infty}(1-q^n)$ is the Dedekind eta function.
The functions $\chi_{1/2}$ resp.~$\chi_8$ are modular functions (with character)
for the theta group $\Gamma_{\theta}=\langle S,T^2\rangle$ respectively for the
full modular group ${\rm PSL}(2,\Z)=\langle S,T\rangle$, where the matrices $S$ 
and $T$ correspond to the modular transformations $\tau\mapsto -\frac{1}{\tau}$ 
and $\tau\mapsto\tau+1$, i.e. $\chi_V$ is a modular function for 
$\Gamma_{\theta}$. The modular curve $\overline{\H/\Gamma_{\theta}}$ has two 
cusps represented by $\infty$ and~$1$. Mainly by definition the character 
$\chi_{V'}$ of the shadow $V'$ is given by the expansion of the modular function
$\chi_V$ in the cusp~$1$ (see~\cite{Ho-dr}, Def.~2.2.4):
\begin{equation}\label{chishadow}
\chi_{V'}(\tau)=\alpha\,e^{2\pi i\frac{c}{24}}\chi_V(-\frac{1}{\tau}
+1),
\end{equation}
with $\alpha=1$ if $c\in\Z$ and $\alpha=\frac{1}{\sqrt{2}}$ if $c\in\
\Z+\frac{1}{2}$. Since
\begin{equation}\label{chi12}
\chi_{\VF'}=\frac{1}{\sqrt{2}}e^{\frac{2\pi i}{48}}\chi_{\VF}(-\frac{1}{\tau}+1)
=q^{\frac{1}{24}}(1+q+q^2+2\,q^3 +\cdots)
\end{equation} 
and
\begin{equation}\label{chi8}
e^{\frac{2\pi i}{3}}\chi_8(-\frac{1}{\tau}+1)=\chi_8(\tau)=q^{-\frac{1}{3}}
(1+248\,q+4124\,q^2+\cdots)
\end{equation}
it follows that the character of the shadow can be written as
\begin{equation}\label{chishadowpol}
\chi_{V'}=\alpha\, P_V(\sqrt{2}\,\chi_{\VF'},\chi_8)={\rm dim\,}V_{h(V')}\cdot
q^{-\frac{c}{24}+h(V')}+\cdots\,.
\end{equation}
From this and the expansions~(\ref{chi12}) and (\ref{chi8}) we conclude that 
for the coefficients of the polynomial 
$P_V(x,y)=\sum_{i=0}^{\left[\frac{c}{8}\right]} A_i\cdot x^{2c-16i}y^i$ one gets
\begin{equation}\label{zero}
A_k=0 \quad \hbox{for $k>\frac{c}{8}-h(V')$.}
\end{equation}
Since the series $\chi_V=P_V(\chi_{1/2},\chi_8)$ starts with $q^{-\frac{c}{24}}$
we get the relation $\sum_{i=0}^{\left[\frac{c}{8}\right]}A_i=1$.
The assumption $h(V')\geq \frac{c}{8}$ implies $P_V(x,y)=x^{2c}$
and one has $h(V')=\frac{c}{8}$. Finally $\chi_V=P_V(\chi_{1/2},\chi_8)=
\chi_{1/2}^{2c}$ gives ${\rm dim\,}V_{1/2}=2c$ and the SVOA $V$ is completely
generated by $V_{1/2}$, i.e.~$V\cong\VF^{\otimes 2c}$.\qed

Similar as in the case of codes and lattices every unitary SVOA $V$ can
be decomposed as a tensor product $W\otimes\VF^{\otimes r}$, where 
$\VF^{\otimes r}$ is the SVOA generated by the elements of conformal weight 
$\halb$ in $V$ and $W$ is a SVOA with no nonzero elements of weight $\halb$.
If $V$ is self-dual, so is $W$ (see~\cite{Ho-dr}, Th.~2.2.8). 
For the shadow of a tensor product we have the the following result:
\begin{lemma}\nonumber Let $V$ and $W$ self-dual unitary ``very nice'' rational SVOAs.
Then $V\otimes W$ is ``very nice'' and for the shadows one has 
$(V\otimes W)'= m\,V'\otimes W'$ with $m=2$ if the ranks of $V$ and $W$
are both half-integral and $m=1$ else.
\end{lemma}
{\bf Proof:} Only the property ``very nice'' has not been completely
shown in~\cite{Ho-dr}, especially it remaines to prove that the even sub-VOA
$$(V \otimes W)_{(0)}=V_{(0)} \otimes W_{(0)}\oplus V_{(1)} \otimes W_{(1)}$$
has the right ``braided fusion algebra'' and gives the right 
$\slz$-representation. Let \hbox{$c$, $d\in \halb\Z$} be the ranks of $V$ and 
$W$. As discussed at the end of the introduction, the ``braided fusion algebras''
of $V_{(0)}$ and $W_{(0)}$ are the same as for $(\VF^{\otimes 2c})_0$ and 
$(\VF^{\otimes 2d})_0$; the ``braided fusion algebra'' of 
$V_{(0)} \otimes W_{(0)}\oplus V_{(1)} \otimes W_{(1)}$ is uniquely determined 
since $V\otimes W$ is unique as $V_{(0)} \otimes W_{(0)}$-module.
So the structure of the ``braided fusion algebras'' is the same as for 
$(\VF^{\otimes 2c}\otimes\VF^{\otimes 2d})_{(0)}=(\VF^{\otimes 2(c+d)})_{(0)}$, 
but $\VF^{\otimes 2(c+d)}$ is known to be ``very nice''.
In particular one gets the stated relations between the shadows and
the structure of the $\slz$-representation follows also. 
\qed

Applying the lemma to $V=W\otimes \VF^{\otimes r}$ we get $h(V')=
h(W')+\frac{r}{16}$, i.e.~the difference $\frac{1}{8}{\rm rk\,}(W)-h(W')$
doesn't change if we tensorize $W$ with copies of $\VF$ and so we can restrict 
us in the following theorem to self-dual SVOAs $V$ with ${\rm dim\,}V_{1/2}=0$.   
\begin{satz}\label{thlong}
Let $V$ be a self-dual unitary ``very nice'' rational SVOA $V$ of rank $c$ 
without nonzero elements of weight $\halb$. Then one has the following:
\begin{description}{}{}
\item i) The dimension of the weight $1$ part $V_1$ is at least $2c\,(23\halb-c)$.
\item ii) The equality holds if and only if the minimal weight of 
the shadow $V'$ is equal to $\frac{c}{8}-1$.
\item iii) In this case the number of linear independent vectors in the shadow
of weight $\frac{c}{8}-1$ is $2^{[c]-11}c$.
\end{description}
\end{satz}
{\bf Proof:} As in the proof of Theorem~\ref{thVF} we use the the polynomial 
$P_V(x,y)$. Assume first that the minimal weight of $V'$ is greater or
equal then $\frac{c}{8}-1$. Equation~(\ref{zero}) shows that $P_V(x,y)$ is a 
linear combination of $x^{2c}$ and $x^{2c-16}y$. From the 
expansions
$$\chi_{1/2}=q^{-\frac{1}{48}}(1+q^{\frac{1}{2}}+q+q^{\frac{3}{2}}+q^2+\cdots)$$
 and~(\ref{chi8}) we find
\begin{eqnarray}\label{shadowlong}
\chi_V& =& \chi_{1/2}^{2c}-\frac{c}{8}\,\chi_{1/2}^{2c-16}
\left(\chi_{1/2}^{16}-\chi_8\right) \\
& =& q^{-\frac{c}{24}}\left(1+0\cdot q^{\frac{1}{2}}+2c\,(23\halb-c)q +\cdots
\right). \nonumber
\end{eqnarray}
Thus ${\rm dim\,} V_1=2c\,(23\halb-c)$ and we have proved one direction of ii).
For the converse we can us restrict to $c<24$. In this case $P_V(x,y)$ hast at
most $3$ terms and is completely determined by the first $3$ coefficients of
$\chi_V$:
$$\chi_V=\chi_{1/2}^{2c}-\frac{c}{8}\,\chi_{1/2}^{2c-16}\left(\chi_{1/2}^{16}-
\chi_8\right)+ \frac{{\rm dim\,}V_1-2c\,(23\halb-c)}{16^2}
\chi_{1/2}^{2c-32}\left(\chi_{1/2}^{16}-\chi_8\right)^2. $$
Equations~(\ref{chi12}), (\ref{chi8}) and (\ref{chishadowpol}) give us then 
for the shadow
$$ {\rm dim\,} V'_{\frac{c}{8}-2}=\alpha\,2^{c-24}\left({\rm dim\,}V_1-
2c\,(23\halb-c)\right).$$
Since ${\rm dim\,}V'_{\frac{c}{8}-2}\geq0$ this implies ${\rm dim\,} V_1\geq
2c\,(23\halb-c)$ and we have proved part i) of the theorem. The equality
holds if an only if ${\rm dim\,}V'_{\frac{c}{8}-2}$ vanishes, which gives
the reverse implication of ii).
Part iii) is finally proved by using~(\ref{chishadowpol}) and (\ref{shadowlong})
to obtain
\begin{eqnarray*}
\chi_{V'}/\alpha & = & (\sqrt{2}\,\chi_{\VF'})^{2c}-\frac{c}{8}
(\sqrt{2}\,\chi_{\VF'})^{2c-16}\left(
(\sqrt{2}\,\chi_{\VF'})^{16}-\chi_8\right) \\
& = & q^{-\frac{c}{24}}\left(\frac{c}{8}\cdot 2^{c-8}\cdot q^{\frac{c}{8}-1}+
\cdots\right),
\end{eqnarray*}
so indeed ${\rm dim\,}V'_{\frac{c}{8}-1}=2^{[c]-11}\cdot c$ which had to be 
shown. \qed

\medskip
To find all the SVOAs $V$ as in Theorem~\ref{thlong} with $h(V')=
\frac{c}{8}-1$ one needs the table of all self-dual SVOAs $V$ of rank $c$
smaller than $24$. By the methods introduced in~\cite{Ho-dr}, Chapter~3 one can 
compute the table and this work is in preparation~\cite{Hoehn2}. For $c<16$ 
every $V$ with ${\rm dim\,}V_{1/2}=0$ satisfies $h(V')=\frac{c}{8}-1$ and the 
list of this SVOAs (making some additional assumptions) together with the
characters of $V$ and $V'$ were already given in~\cite{Ho-dr} 
(Theorem~3.2.4 and Table~5.3); confirming the the values for ${\rm dim\,} V_1$
and ${\rm dim\,}V'_{\frac{c}{8}-1}$ as stated in Theorem~\ref{thlong}. 
Using also~\cite{Hoehn2} we get the following list:

$$
\begin{array}{c|cccccccccccc}
c   & 8 & 12 & 14 & 15 & 15\halb & 16 &   17 & 17\halb & 18 & 18\halb \\ \hline
{\rm dim\,} V_1 & 248 & 276 & 266 & 255 & 248 & 240  & 221 & 210 & 198 & 185 \\ 
\hline
\widetilde{V}_1 & E_8 & D_{12} & E_7^2 & A_{15} & E_{8,2} & D_8^2 &
A_{11}E_6 & C_{10} & D_6^3 & E_{7,2}F_4   \\ 
\end{array}
$$ $$
\begin{array}{c|cccccccccccc}
c   & 19 & 19\halb & 20 & 20\halb & 21 & 21\halb & 22 &
22\halb & 23 & 23\halb \\ \hline 
{\rm dim\,} V_1  & 171 & 156 & 140 & 123 & 105 & 86 & 66 & 45 & 23 & 0 \\ \hline
\widetilde{V}_1 &A_7^2D_5 & D_{8,2}B_4 & D_4^5 & 
A_{9,2}A_4 & A_3^7 & D_{4,2}^2 C_2^3 & A_1^{22} & A_{1,2}^{15} & U_1^{23} & 
0\subset \VB  \\ 
\end{array}
$$

The first row gives the rank, the second the dimension of $V_1$ and the last
row contains one example of a SVOA meeting the bound of 
Theorem~\ref{thlong},~i) labelled by the sub-VOA $\widetilde{V}_1$ generated 
by $V_1$, this is the VOA associated to an integrable highest weight 
representation of the affine Lie algebra coming from the reductive Lie 
algebra $V_1$. The first index of a component denotes the rank of the simple
Lie subalgebra, the second (if it exist) the level of the highest weight 
representation. This index is omitted if the level is $1$.

As it is the case for lattices and codes the following holds for SVOAs: 
If there exist for a $c$ in the range  $16\leq c<24$ a SVOA of rank $c$ 
with ${\rm dim\,}V_{1/2}=0$ at all,
then there exist such a SVOA meeting the bound 
${\rm dim\,} V_1=2c\,(23\halb-c)$. 
The first SVOA $V_{E_8}$ in the list is even a VOA. 
Especially interesting is the last one. The SVOA
$\VB$ of rank $23\frac{1}{2}$ is the {\it shorter moonshine module\/} or 
{Baby Monster SVOA\/} and was constructed in~\cite{Ho-dr}, Chapter~4. On $\VB$
acts the group $2 \times B$ by automorphisms, where $B$ is the Baby Monster,
the second largest sporadic simple group. It is conjectured that $\VB$ is
the unique self-dual unitary ``very nice'' rational SVOA with ${\rm dim\,}
V_{1/2}={\rm dim\,} V_1=0$.

\section{Relations to codes and lattices}

The analogs of Theorem~\ref{thVF} and~\ref{thlong} and of the above table for
self-dual codes and lattices as given in~\cite{El-Z,El-shadow} can easily be deduced 
from the results for SVOAs. Since the reduction from codes to lattices was 
described in~\cite{El-shadow} we show only the (analog) reduction from SVOAs
to lattices. For the notation see also~\cite{El-shadow}.
\begin{korollar}
Let $L$ be a self-dual lattice in $\R^n$. Then
the analogs of Theorem~\ref{thVF} and~\ref{thlong} hold for $L$:
\begin{description}
\item
(1) The shortest characteristic vector of $L$ has norm $n$ if and only if
$L\cong \Z^n$.
\item 
(2) Assume that $L$ has no vectors of norm $1$. Then:
\begin{description}{}{}
\item i) $L$ has at least $2n(23-n)$ vectors of norm $2$.
\item ii) Equality holds if and only if $L$ has no characteristic vectors of 
norm smaller then $n-8$.
\item iii) In that case the number of characteristic vectors of length 
exactly $n-8$ is $2^{n-11}n$.
\end{description}
\end{description}
\end{korollar}
{\bf Proof:} To every integral lattice $L$ of rank $n$ one can construct a
SVOA $V_L$ of rank $c=n$ whose irreducible modules $M_{\lambda}$ are indexed
be the cosets $\lambda\in L^*/L$ of the dual lattice modulo the lattice.
The character of $M_{\lambda}$ is determined by the theta series of $L+\lambda$:
\begin{equation}\label{chiVL}
\chi_{M_{\lambda}}= \big(\sum_{x\in L+\lambda} q^{\frac{1}{2}\langle x,x
\rangle}\big)\cdot \big(q^{1/24}\prod_{l=1}^{\infty}(1-q^l)\big)^{-n}.  
\end{equation}
For self-dual lattices one has $V_L'=M_{L'}$ where
$L'=L^*\setminus L$ is the shadow of $L$ and (\ref{chiVL}) shows that
$h(V_L')$ is equal to $\frac{1}{8}$ times the minimal norm of a characteristic
vector of $L$.
Since $V_{\Z}\cong\VF^{\otimes 2}$ and $V_{K\oplus L}\cong V_K\otimes V_L$ for
two integral lattices $K$ and $L$ part (1) follows from Theorem~\ref{thVF}.

For (2) we note first, that ${\rm dim\,}(V_L)_{1/2}=0$ if and only if $L$
has no vectors of norm $1$. From~(\ref{chiVL}) we get that
${\rm dim\,}(V_L)_1$ is equal to $n$ plus the numbers of vectors in $L$ of
norm $2$ and
${\rm dim\,}(V_L')_{\frac{c}{8}-1}$ is the number of characteristic vectors 
of $L$ of norm $(n-8)$.
Together with the identity $2n(23\halb-n) = 2n(23-n)+n$ part (2) of the 
corollary follows. \qed

A SVOA $V$ can only be come from a lattice by the above construction if its rank
is integral and the affine Kac-Moody sub-VOA $\widetilde{V}_1$ is a 
tensor product of level-$1$ representations of affine Kac-Moody algebras 
corresponding to the ``simply laced'' Lie algebras $A_n$, $D_n$, $E_6$,
$E_7$ and $E_8$. (One has $\widetilde{(V_L)}_1=V_{L_0}$, where $L_0$ is the 
root sublattice of $L$.) In fact every SVOA of Theorem~\ref{thlong} with 
$h(V')=\frac{c}{8}-1$
and the above properties comes from a unique lattice, thus the corresponding
list for lattices can be read off from the list of such SVOAs. 
 
\medskip
{\it Remarks:\/} One can define Kleinian codes as a natural fourth step before 
binary codes, lattices and vertex operator algebras and then again the analogous
theorems hold~\cite{Ho-kleinian}.

\smallskip
The analog of Theorem~\ref{thVF} for lattices was used 
in~\cite{El-Z} as the final step in the proof of Donaldson famous theorem 
about positive definite intersection forms of closed differentiable simply 
connected $4$-manifolds with the help of the Seiberg-Witten invariants. 
Theorem~\ref{thVF} relates quantum field theory in another way to the new proof
of this theorem.


\bigskip
\noindent
Current Address: {\sc
Universit\"at Freiburg,
Mathematisches Institut,
Eckerstr.~1,
79104 Freiburg,
Germany;} E-mail: {\tt gerald@mathematik.uni-freiburg.de}
 
\end{document}